\documentclass[11pt]{article}

\usepackage{graphicx,psfrag}
\usepackage{psfrag,eepic,epsfig}
\usepackage{sectsty}
\usepackage{setspace}
\usepackage{color}

\usepackage{setspace}
\usepackage{amsmath, epsfig, cite, ifpdf}
\usepackage{amssymb}
\usepackage{amsfonts}
\usepackage{latexsym}
\usepackage{graphicx,psfrag,eepic,rawfonts}
\usepackage{enumerate}
\usepackage{indentfirst}
\input{epsf}

\newtheorem{thm}{Theorem}[section]
\newtheorem{prop}[thm]{Proposition}

\newtheorem{lem}[thm]{Lemma}

\newcommand{\qed}{{\hfill\rule{4pt}{7pt}}}
\def\pf{\noindent {\it Proof.} }

\numberwithin{equation}{section}

\makeatletter \@addtoreset{equation}{section} \makeatother

\title {\bf $4$-Regular oriented graphs with \\optimum skew
energy\footnote{Supported by NSFC and the ``973" program. }}

\author{
{\small Xiaolin Chen, Xueliang Li, Huishu Lian}\\
{\small Center for Combinatorics and LPMC-TJKLC}\\
{\small Nankai University, Tianjin 300071, P.R. China}\\
{\small E-mail: chxlnk@163.com; lxl@nankai.edu.cn; lhs6803@126.com}
   }
\date{}
\begin{document}

\maketitle

\begin{abstract}
Let $G$ be a simple undirected graph, and $G^\sigma$ be an oriented
graph of $G$ with the orientation $\sigma$ and skew-adjacency matrix
$S(G^\sigma)$. The skew energy of the oriented graph $G^\sigma$,
denoted by $\mathcal{E}_S(G^\sigma)$, is defined as the sum of the
absolute values of all the eigenvalues of $S(G^\sigma)$. In this
paper, we characterize the underlying graphs of all $4$-regular
oriented graphs with optimum skew energy and give orientations of
these underlying graphs such that the skew energy of the resultant
oriented graphs indeed attain optimum. It should be pointed out that
there are infinitely many 4-regular connected optimum skew energy
oriented graphs, while the $3$-regular case only has two graphs: 
$K_4$ the complete graph on $4$ vertices and $Q_3$ the hypercube.


\end{abstract}

\section{Introduction}

Let $G$ be a simple undirected graph with vertex set
$V(G)=\{v_1,v_2,\ldots,v_n\}$, and let $G^\sigma$ be an oriented
graph of $G$ with the orientation $\sigma$, which assigns to each
edge of $G$ a direction so that the induced graph $G^\sigma$ becomes
an oriented graph or a directed graph. Then $G$ is called the underlying
graph of $G^\sigma$. The skew-adjacency matrix of $G^\sigma$ is the
$n\times n$ matrix $S(G^\sigma)=[s_{ij}]$, where $s_{ij}=1$ and
$s_{ji}=-1$ if $\langle v_i,v_j\rangle$ is an arc of $G^\sigma$,
otherwise $s_{ij}=s_{ji}=0$. The skew energy \cite{ABC} of $G^\sigma$,
denoted by $\mathcal{E}_S(G^\sigma)$, is defined as the sum of the
absolute values of all the eigenvalues of $S(G^\sigma)$. Obviously,
$S(G^\sigma)$ is a skew-symmetric matrix, and thus all the eigenvalues
are purely imaginary numbers.

In theoretical chemistry, the energy of a given molecular graph is
related to the total $\pi$-electron energy of the molecule
represented by that graph. Consequently, the graph energy has some
specific chemistry interests and has been extensively studied, since
the concept of the energy of simple undirected graphs was introduced
by Gutman in \cite{G}. We refer the survey \cite{GLZ} and the book
\cite{LSG} to the reader for details. Up to now, there are various
generalizations of the graph energy, such as the Laplacian energy,
signless Laplacian energy, incidence energy, distance energy, and
the Laplacian-energy like invariant for undirected graphs, and the
skew energy and skew Laplacian energy for oriented graphs.

Adiga et al. \cite{ABC} first defined the skew energy of an oriented
graph, and investigated some properties of the skew energy. Then,
Shader et al. \cite{Shader} studied the relationship between the
spectra of a graph $G$ and the skew-spectra of an oriented graph
$G^\sigma$ of $G$, which would be helpful to the study of the
relationship between the energy of $G$ and the skew energy of
$G^\sigma$. Hou and Lei \cite{HL} characterized the coefficients of
the characteristic polynomial of the skew-adjacency matrix of an
oriented graph. Moreover, other bounds and extremal graphs of some
classes of oriented graphs have been established. In \cite{HSZ} and
\cite{HSZ2}, Hou et al. determined the oriented unicyclic graphs
with minimal and maximal skew energy and the oriented bicyclic
graphs with minimal and maximal shew energy, respectively. The skew
energy of orientations of hypercubes were discussed by Tian
\cite{Tian}. Later, Gong and Xu \cite{GX} characterized the
3-regular oriented graphs with optimum skew energy. Recently, we
\cite{CLL} studied the skew energy of random oriented graphs.

Back to the paper Adiga et al. \cite{ABC}, where they derived a
sharp upper bound for the skew energy of an oriented graph
$G^\sigma$ in terms of the order $n$ and the maximum degree $\Delta$
of $G^\sigma$, that is,
\begin{equation*}
\mathcal{E}_S(G)\leq n\sqrt{\Delta}\,.
\end{equation*}
They showed that the equality holds if and only if
$S(G^\sigma)^{T}S(G^\sigma)=\Delta I_n$, which implies that
$G^\sigma$ is $\Delta$-regular. In the following, we will call an
oriented graph $G^\sigma$ on $n$ vertices with maximum degree
$\Delta$ an {\it optimum skew energy oriented graph} if
$\mathcal{E}_S(G^\sigma)= n\sqrt{\Delta}$. A natural question is
proposed in \cite{ABC}:

Which $k$-regular graphs on $n$ vertices have orientations
$G^\sigma$ with $\mathcal{E}_S(G^\sigma)= n\sqrt{\Delta}$, or
equivalently, $S(G^\sigma)^{T}S(G^\sigma)=\Delta I_n$ ?

\noindent In the same paper, they answer the question for $k=1$ and
$k=2$. They showed that a $1$-regular graph on $n$ vertices has an
orientation with $S(G^\sigma)^{T}S(G^\sigma)=I_n$ if and only if $n$
is even and it is $\frac{n}{2}$ copies of $K_2$; while a $2$-regular
graph on $n$ vertices has an orientation with
$S(G^\sigma)^{T}S(G^\sigma)=2 I_n$ if and only if $n$ is a multiple
of $4$ and it is a union of $\frac{n}{4}$ copies of $C_4$. Later,
Gong and Xu \cite{GX} characterized all $3$-regular connected
oriented graphs on $n$ vertices with
$S(G^\sigma)^{T}S(G^\sigma)=3I_n$, which in fact are only two
special graphs, the complete graph $K_4$ and the hypercube $Q_3$.

In this paper, we further consider the above question. We
characterize all $4$-regular connected graphs $G$ that have oriented
graphs $G^\sigma$ with $S(G^\sigma)^{T}S(G^\sigma)=4 I_n$. It should
be noted that the $4$-regular case is more complicated than the
$3$-regular case, and in fact, there are infinitely many 4-regular
connected optimum skew energy oriented graphs.

\section{Preliminaries}

In this section, we do some preparations with some notations and a
few known results. Besides, we also get some intuitive conclusions
that will be frequently used in the sequel of the paper.

Let $G=G(V,E)$ be a graph with vertex set $V$ and edge set $E$. For
any $v\in V$, denote by $d_G(v)$ and $N_G(v)$ the degree and
neighborhood of $v$ in $G$, respectively. For any subset $S\subseteq
V$, $G[S]$ denotes the subgraph of $G$ induced by $S$. For a given
orientation $\sigma$ of $G$, the resultant oriented graph is denoted
by $G^\sigma=(V(G^\sigma),\Gamma(G^\sigma))$ and the skew-adjacency
matrix of $G^\sigma$ by $S(G^\sigma)$.

The following result is due to Adiga et al. \cite{ABC}.
\begin{lem}\cite{ABC}\label{neigh}
Let $S(G^\sigma)$ be the skew-adjacency matrix of an oriented
graph $G^\sigma$. If $S(G^\sigma)^TS(G^\sigma)=kI$, then $|N(u)\cap N(v)|$
is even for any two distinct vertices $u$ and $v$ of $G^\sigma$.
\end{lem}

Since our paper focuses on the investigation of $4$-regular graphs,
the following result is more directly applied, which is in fact
implied in Lemma \ref{neigh}.
\begin{prop}\label{4neigh}
Let $G^\sigma$ be a $4$-regular oriented graph with skew-adjacency
matrix $S(G^\sigma)$. If $S(G^\sigma)^TS(G^\sigma)=4I$, then the
underlying graph $G$ satisfies that $|N(u)\cap N(v)|\in \{0,2\}$ for
any two adjacent vertices $u$ and $v$ and $|N(u)\cap N(v)|\in
\{0,2,4\}$ for any two non-adjacent vertices $u$ and $v$.
\end{prop}

Let $W=u_1u_2\cdots u_k$ (perhaps $u_i=u_j$ for $i\neq j$) be a walk
from $u_1$ to $u_k$ and $\widehat{W}$ be the inverse walk of $W$
obtained from $W$ by replacing the ordering of vertices by its
inverses, i.e., $\widehat{W}=u_ku_{k-1}\cdots u_1$. The sign of $W$
is defined as
\begin{equation*}
\text{sgn}(W)=\prod_{i=1}^{k-1} s_{u_iu_{i+1}}.
\end{equation*}
It is easy to check that
\begin{equation*}
\text{sgn}(\widehat{W})=
\begin{cases}
\text{sgn}(W) & \text{if\,\,} l(W) \text{\,\,is even},\\
-\text{sgn}(W) & \text{if\,\,} l(W) \text{\,\,is odd},
\end{cases}
\end{equation*}
where $l(W)$ denotes the length of the walk $W$. Moreover, let
$w_{uv}^+(k)$ and $w_{uv}^-(k)$ denote the number of all positive
walks and negative walks starting from $u$ and ending at $v$ with
length $k$, respectively.

Gong and Xu \cite{GX} obtained the following result on the
relationship between the entries of $S^k$ and the number of walks
between any pair of ordered vertices.
\begin{lem}\cite{GX}\label{walk}
Let $S$ be the skew-adjacency matrix of an oriented graph $G^\sigma$
and $(u,v)$ be an arbitrary pair of ordered vertices of $G^\sigma$.
Then
\begin{equation*}
(S^k)_{uv}=w_{uv}^+(k)-w_{uv}^-(k)
\end{equation*}
holds for any positive integer $k$.
\end{lem}

For regular graphs, the following proposition is immediate.
\begin{prop}\label{2-walk}
Let $G^\sigma$ be a $k$-regular oriented graph with skew-adjacency
matrix $S$. Then $S^TS=kI$ if and only if for any two distinct
vertices $u$ and $v$ of $G^\sigma$,
\begin{equation*}
w_{uv}^+(2)=w_{uv}^-(2).
\end{equation*}
\end{prop}

Throughout this paper, we just need to consider connected graphs
and connected oriented graphs due to the following basic lemma.
Recall that the union $G_1^\sigma\cup G_2^\sigma$ of two disjoint
oriented graphs $G_1^\sigma=(V_1,\Gamma_1)$ and $G_2^\sigma=(V_2,\Gamma_2)$
is the oriented graph $G^\sigma=(V,\Gamma)$ where $V=V_1\cup V_2$
and $\Gamma=\Gamma_1\cup \Gamma_2$.
\begin{lem}\cite{Tian}\label{connected}
Let $G_1^\sigma$, $G_2^\sigma$ be two disjoint oriented graphs of
order $n_1$, $n_2$ with skew-adjacency matrices $S(G_1^\sigma)$,
$S(G_2^\sigma)$, respectively. Then for some positive integer $k$,
$S(G_1^\sigma)^TS(G_1^\sigma)=kI_{n_1}$ and $S(G_2^\sigma)^TS(G_2^\sigma)=kI_{n_2}$
if and only if the skew-adjacency matrix $S(G_1^\sigma \cup G_2^\sigma)$
of the union $G_1^\sigma\cup G_2^\sigma$ satisfies
$S(G_1^\sigma\cup G_2^\sigma)^TS(G_1^\sigma\cup G_2^\sigma)=kI_{n_1+n_2}$.
\end{lem}

\begin{figure}[h,t,b,p]
\begin{center}
\scalebox{1}[1]{\includegraphics{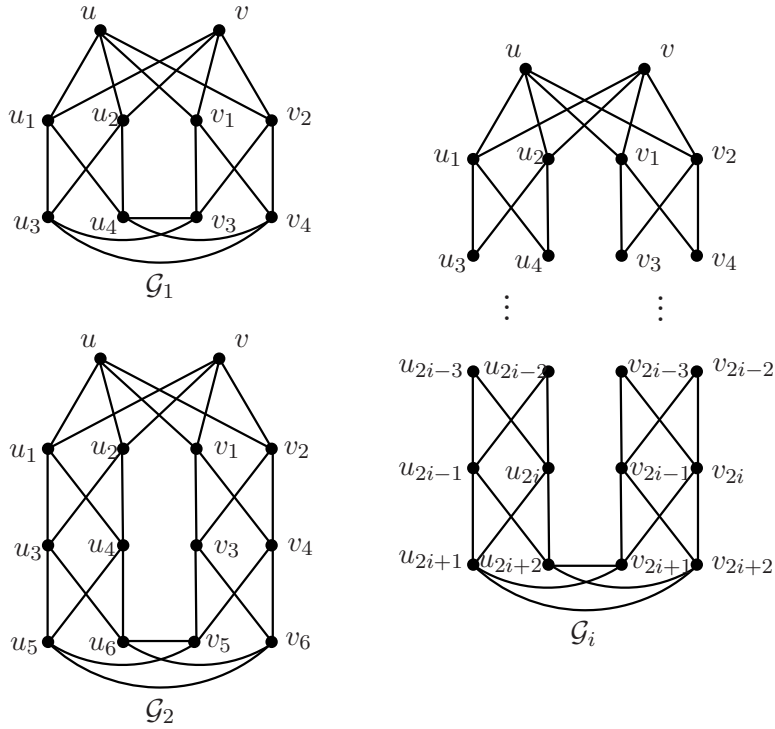}}
\end{center}
\caption{The graph class $\mathcal{G}_i$ for any  positive integer $i$}\label{Fig1}
\end{figure}

\begin{figure}[h,t,b,p]
\begin{center}
\scalebox{1}[1]{\includegraphics{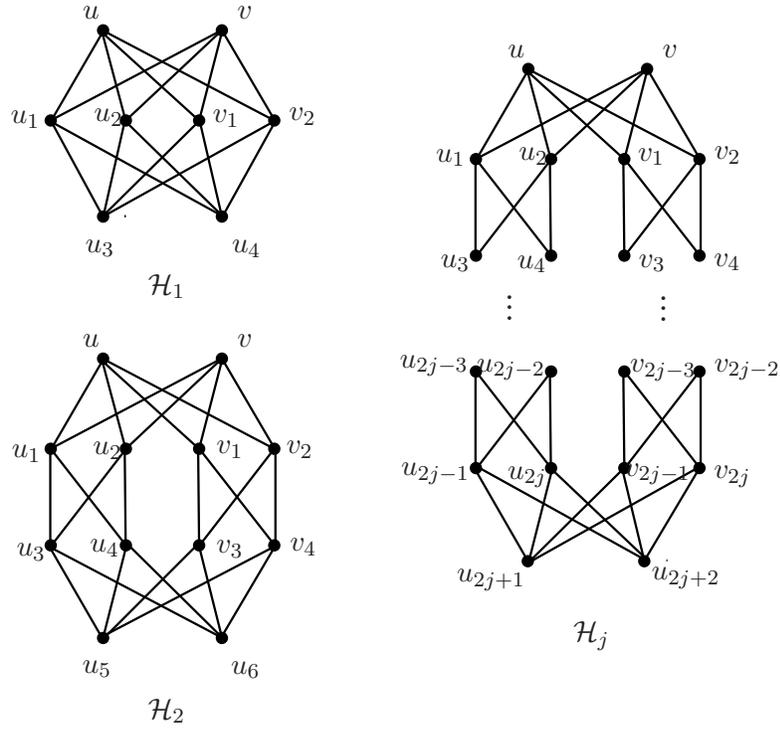}}
\end{center}
\caption{The graph class $\mathcal{H}_j$ for any  positive integer $j$}\label{Fig2}
\end{figure}

We end this section by recursively defining two graph classes $\mathcal{G}_i$
and $\mathcal{H}_j$ for all positive integers $i$ and $j$, depicted
in Figure \ref{Fig1} and Figure \ref{Fig2}, respectively.

For the graph class $\mathcal{G}_i$, we define the initial graph
$\mathcal{G}_1=(V(\mathcal{G}_1),E(\mathcal{G}_1))$, where
\begin{align*}
V(\mathcal{G}_1)=&\{u,v\}\cup\{u_1,u_2,v_1,v_2\}\cup \{u_3,u_4,v_3,v_4\},\\
E(\mathcal{G}_1)=&\{(u,u_1),(u,u_2),(u,v_1),(u,v_2),(v,u_1),(v,u_2),(v,v_1),(v,v_2)\}\\
&\cup \{(u_1,u_3),(u_1,u_4),(u_2,u_3),(u_2,u_4),(v_1,v_3),(v_1,v_4),(v_2,v_3),(v_2,v_4)\}\\
&\cup \{(u_3,v_4),(v_4,u_4),(u_4,v_3),(v_3,u_3)\}.
\end{align*}
Suppose that $\mathcal{G}_{i-1}$ is well defined. Below we will give
the definition of
$\mathcal{G}_i=(V(\mathcal{G}_i),E(\mathcal{G}_i))$.
\begin{align*}
V(\mathcal{G}_i)=&V(\mathcal{G}_{i-1})\cup\{u_{2i+1},u_{2j+2},v_{2i+1},v_{2i+2}\},\\
E(\mathcal{G}_i)=&E(\mathcal{G}_{i-1})\setminus \{(u_{2i-1},v_{2i}),(v_{2i},u_{2i}),(u_{2i},v_{2i-1}),(v_{2i-1},u_{2i-1})\}\\
&\cup \{(u_{2i-1},u_{2i+1}),(u_{2i-1},u_{2i+2}),(u_{2i},u_{2i+1}),(u_{2i},u_{2i+2})\}\\
&\cup \{(v_{2i-1},v_{2i+1}),(v_{2i-1},v_{2i+2}),(v_{2i},v_{2i+1}),(v_{2i},v_{2i+2})\}\\
&\cup \{(u_{2i+1},v_{2i+2}),(v_{2i+2},u_{2i+2}),(u_{2i+2},v_{2i+1}),(v_{2i+1},u_{2i+1})\}.
\end{align*}
Observe that $|V(\mathcal{G}_i)|=4i+6$.

For the other graph class $\mathcal{H}_j$, the initial graph $\mathcal{H}_1$
is defined as $\mathcal{H}_1=(V(\mathcal{H}_1),E(\mathcal{H}_1))$, where
\begin{align*}
V(\mathcal{H}_1)=&\{u,v\}\cup \{u_1,u_2,v_1,v_2\}\cup \{u_3,u_4\},\\
E(\mathcal{H}_1)=&\{(u,u_1),(u,u_2),(u,v_1),(u,v_2),(v,u_1),(v,u_2),(v,v_1),(v,v_2)\}\\
&\cup \{(u_1,u_3),(u_1,u_4),(u_2,u_3),(u_2,u_4),(v_1,u_3),(v_1,u_4),(v_2,u_3),(v_2,u_4)\}.
\end{align*}
Suppose now that we have given the definition of
$\mathcal{H}_{j-1}$. Then
$\mathcal{H}_j=(V(\mathcal{H}_j),E(\mathcal{H}_j))$ is defined as
follows.

\begin{align*}
V(\mathcal{H}_j)=&V(\mathcal{H}_{j-1})\cup\{v_{2j-1},v_{2j},u_{2j+1},u_{2j+2}\},\\
E(\mathcal{H}_j)=&E(\mathcal{H}_{j-1})\setminus \{(v_{2j-3},u_{2j-1}),
(v_{2j-3},u_{2j}),(v_{2j-2},u_{2j-1}),(v_{2j-2},u_{2j})\}\\
&\cup \{(v_{2j-3},v_{2j-1}),(v_{2j-3},v_{2j}),(v_{2j-2},v_{2j-1}),(v_{2j-2},v_{2j})\}\\
&\cup \{(u_{2j-1},u_{2j+1}),(u_{2j-1},u_{2j+2}),(u_{2j},u_{2j+1}),(u_{2j},u_{2j+2})\}\\
&\cup \{(v_{2j-1},u_{2j+1}),(v_{2j-1},u_{2j+2}),(v_{2j},u_{2j+1}),(v_{2j},u_{2j+2})\}.
\end{align*}
Obviously, $|V(\mathcal{H}_j)|=4j+4$.

\section{Main results}

\begin{figure}[h,t,b,p]
\begin{center}
\scalebox{1}[1]{\includegraphics{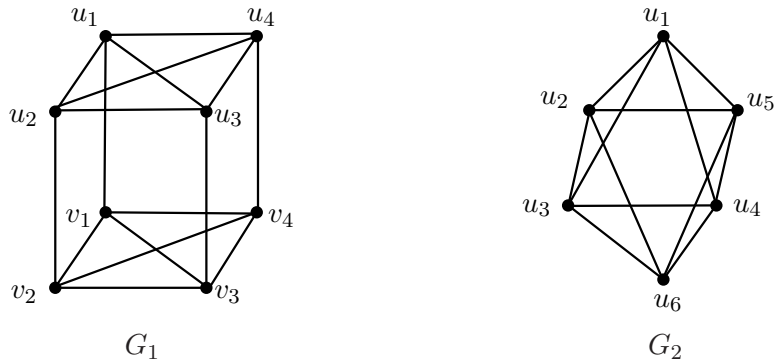}}
\end{center}
\caption{The underlying graphs containing triangles}\label{Fig3}
\end{figure}

In this section, we first characterize the underlying graphs of all
$4$-regular oriented graphs with optimum skew energy. Then we give
orientations of these underlying graphs such that the resultant
oriented graphs have optimum shew energy.
\begin{thm}\label{underlying1}
Let $G^\sigma$ be a $4$-regular oriented graph with optimum skew energy.
If the underlying graph $G$ contains triangles, then $G$ is either $G_1$
or $G_2$ depicted in Figure \ref{Fig3}.
\end{thm}

\pf Let $u_1u_2u_3u_1$ be a triangle in $G$. Since $u_2 \in
N(u_1)\cap N(u_3)$, there is another common neighbor between $u_1$
and $u_3$ from Proposition \ref{4neigh}, denoted by $u_4$. Observe
that $u_3 \in N(u_1)\cap N(u_2)$. Then by Proposition \ref{4neigh}
again, there is another vertex in $N(u_1)\cap N(u_2)$, which is
either $u_4$ or a new vertex, say $u_5$.

Firstly, assume that $u_4\in N(u_1)\cap N(u_2)$, that is,
$(u_2,u_4)\in G$. As $G$ is $4$-regular, $u_1$ has the fourth
neighbor, denoted by $v_1$. We claim that $(v_1,u_2)\notin G$;
otherwise $N(u_1)\cap N(u_2)=\{u_3,u_4,v_1\}$ which contradicts
Proposition \ref{4neigh}. Similarly, we have $(v_1,u_3)\notin G$ and
$(v_1,u_4)\notin G$. We can further obtain that the new vertices
$v_2,v_3$ and $v_4$ are the forth neighbors of $u_2,u_3$ and $u_4$,
respectively, and $(v_i,u_j)\notin G$ for $1\leq i\neq j\leq 4$.
Then we consider $N(v_1)\cap N(u_2)$. Note that $u_1\in N(u_2)\cap
N(v_1)$, $(v_1,u_3)\notin G$ and $(v_1,u_4)\notin G$ by the
discussion above, which forces that $v_2$ becomes another common
neighbor between $v_1$ and $u_2$, i.e., $(v_1,v_2)\in G$. By similar
discussions on $N(v_1)\cap N(u_4)$, $N(v_3)\cap N(u_2)$ and
$N(v_3)\cap N(u_4)$, respectively, we can deduce that $(v_1,v_4)\in
G$, $(v_2,v_3)\in G$ and $(v_3,v_4)\in G$. Noticing that $u_1 \in
N(v_1)\cap N(u_3)$, another common vertex must be $v_3$, since
$d(u_3)=4$ and the degrees of other neighbors of $u_3$ other than
$v_3$ are equal to $4$. By considering $N(u_2)\cap N(v_4)$
similarly, we have $(v_2,v_4)\in G$. Up to now, the degrees of all
vertices of $G$ attain $4$. Hence the underlying graph $G$ is the
graph $G_1$ given in Figure \ref{Fig3}.

\begin{figure}[h,t,b,p]
\begin{center}
\scalebox{1}[1]{\includegraphics{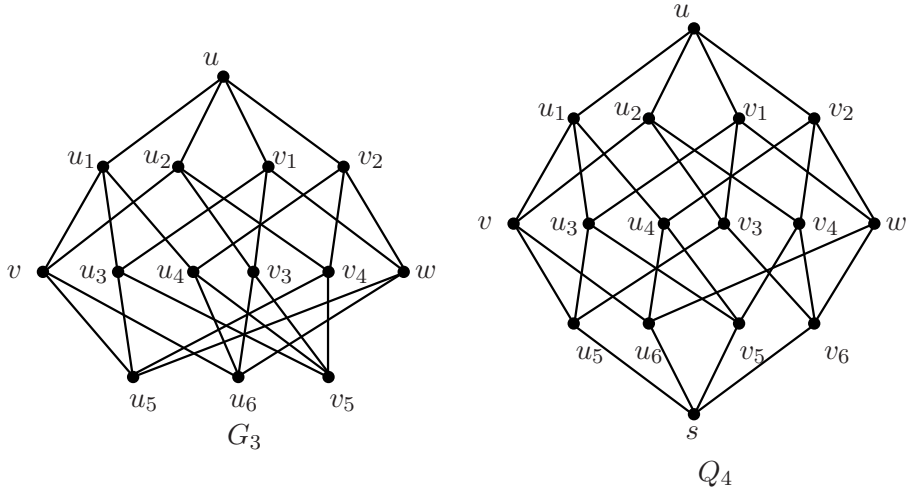}}
\end{center}
\caption{The underlying graphs containing no triangles}\label{Fig4}
\end{figure}

Now we suppose that $N(u_1)\cap N(u_2)$ contains a new vertex $u_5$.
We claim that $(u_2,u_4)\notin G$ and $(u_3,u_5)\notin G$;
otherwise, $N(u_1)\cap N(u_2)=\{u_3,u_4,u_5\}$ or $N(u_1)\cap
N(u_3)=\{u_2,u_4,u_5\}$, a contradiction to Proposition
\ref{4neigh}. Notice that $u_3\in N(u_1)\cap N(u_4)$, $d(u_1)=4$ and
$(u_2,u_4)\notin G$, which implies $(u_4, u_5)\in G$. Since
$d(u_5)=3$ and $(u_3,u_5)\notin G$, $u_5$ has the forth neighbor
$u_6$. Now we consider $N(u_2)\cap N(u_5)$. Combining the
observation that $u_1\in N(u_2)\cap N(u_5)$ with the fact
$(u_2,u_4)\notin G$, we deduce that $u_6\in N(u_2)\cap N(u_5)$. Then
by a similar way, we successively discuss $N(u_2)\cap N(u_3)$ and
$N(u_3)\cap N(u_4)$ and obtain $(u_3,u_6)\in G$ and $(u_4,u_6)\in
G$. It is easy to check that the graph has already been $4$-regular
and is just the graph $G_2$ depicted in Figure \ref{Fig3}.\qed
\begin{thm}\label{underlying2}
Let $G^\sigma$ be a $4$-regular oriented graph with optimum skew
energy. If the underlying graph $G$ is triangle-free, then $G$ is
one of the following graphs: the hypercube $Q_4$ of dimension $4$,
the graph $G_3$, a graph in $\mathcal{G}_i$, or a graph in
$\mathcal{H}_j$; see Figures \ref{Fig1}, \ref{Fig2} and \ref{Fig4}.
\end{thm}

\pf Let $u_1,u_2,v_1$ and $v_2$ be all neighbors of a vertex $u$ in
$G$. Then the induced subgraph $G[\{u_1,u_2,v_1,v_2\}]$ contains no
edge, since the graph $G$ is triangle-free. Denote by $v,u_3,u_4$ be
another three neighbors of $u_1$ other than $u$. Note that $u_1 \in
N(u)\cap N(v)$. By Proposition \ref{4neigh}, there is another one or
three common neighbors in $\{u_2,v_1,v_2\}$ between $u$ and $v$. We
can obtain the same results by considering $N(u)\cap N(u_3)$ and
$N(u)\cap N(u_4)$. Assume that $a_1,a_2$ and $a_3$ are the numbers
of the common neighbors in $\{u_2,v_1,v_2\}$ between $u$ and $v$,
$u$ and $u_3$, $u$ and $u_4$, respectively. Obviously,
$a_1,a_2,a_3\in \{1,3\}$. Without loss of generality, suppose
$a_1\geq a_2\geq a_3$. We discuss the following four cases according
to the values of $a_1,a_2$ and $a_3$.

\textbf{Case 1.} $(a_1,a_2,a_3)=(1,1,1)$.

Without loss of generality, let $(u_2,v)\in G$. Then $(v_1,v)\notin
G$ and $(v_2,v)\notin G$ as $a_1=1$. Observe that $u\in N(u_1)\cap
N(v_1)$, which implies that there is another common neighbor in
$\{u_3,u_4\}$ between $u_1$ and $v_1$. Let $u_3\in N(u_1)\cap
N(v_1)$. Then $(u_2,u_3)\notin G$ and $(v_2,u_3)\notin G$ as
$a_2=1$. By considering $N(u_1)\cap N(v_2)$, we deduce that
$(v_2,u_4)\in G$, since $(v_2,v)\notin G$ and $(v_2,u_3)\notin G$ by
the discussion above. Obviously, $(u_2,u_4)\notin G$ and
$(v_2,u_4)\notin G$ as $a_3=1$. Then it is known that $u_2$ contains
another two neighbors, say $v_3$ and $v_4$. Since $u\in N(u_2)\cap
N(v_1)$ and $(v_1,v)\notin G$, it follows that there is another
common neighbor in $\{v_3,v_4\}$ between $u_2$ and $v_1$. Without
loss of generality, $v_3\in N(u_2)\cap N(v_1)$. Then
$(v_1,v_4)\notin G$; otherwise, $|N(u_2)\cap N(v_1)|=3$ and no other
vertex can be chosen as the forth common neighbor, which is a
contradiction. In view of the observation that $u\in N(u_2)\cap
N(v_2)$ and $(v_2,v)\notin G$, we have that another common neighbor
between $u_2$ and $v_2$ belongs to $\{v_3,v_4\}$. We claim that
$v_3\notin N(u_2)\cap N(v_2)$; otherwise, $N(u)\cap
N(v_3)=\{u_2,v_1,v_2\}$ and there is no other vertex in $N(u)\cap
N(v_3)$, which contradicts Proposition \ref{4neigh}. Therefore,
$v_4\in N(u_2)\cap N(v_2)$. We proceed to have $w$ as the forth
neighbor of $v_1$. By considering $N(v_1)\cap N(v_2)$, we obtain
$(v_2,w)\in G$.

Up to now, we have $d(u)=d(u_1)=d(u_2)=d(v_1)=d(v_2)=4$ and
$d(v)=d(u_3)=d(u_4)=d(v_3)=d(v_4)=d(w)=3$. We claim that the deduced
subgraph $G[\{v,u_3,u_4,v_3,v_4,w\}]$ is empty. Otherwise, the
possible edges are $(v,w),(u_3,v_4)$ and $(u_4,v_3)$ since $G$ is
triangle-free. If $(v,w)\in G$, then $|N(u_1)\cap N(w)|=1$, which is
a contradiction. We thus have $(v,w)\notin G$. Similarly,
$(u_3,v_4)\notin G$ and $(u_4,v_3)\notin G$.

Suppose now that $u_5$ and $u_6$ are the other two neighbors of $v$.
Note that $u_1\in N(v)\cap N(u_3)$. Then we have either
$(u_3,u_5)\in G$ or $(u_3,u_6)\in G$. Without loss of generality,
$(u_3,u_5)\in G$, and hence $(u_3,u_6)\notin G$. Moreover,
$(u_4,u_5)\notin G$, otherwise, $N(u_1)\cap N(u_5)=\{v,u_3,u_4\}$, a
contradiction. By considering $N(u_1)\cap N(u_6)$, we get
$(u_4,u_6)\in G$. Assume that $v_5$ is the forth neighbor of $u_3$.
It is obvious that $u_3\in N(u_1)\cap N(v_5)$, which forces that
$u_4$ becomes another common vertex between $u_1$ and $v_5$. We see
that $v\in N(u_2)\cap N(u_5)$, which indicates that there is another
common neighbor between $u_2$ and $u_5$. It means either
$(v_3,u_5)\in G$ or $(v_4,u_5)\in G$. We discuss the two cases
separately.

On the one hand, if $(v_3,u_5)\in G$, then $(v_4,u_5)\notin G$. It
follows that $(v_4,u_6)\in G$ by considering $N(v)\cap N(v_4)$. We
claim that $(v_3,u_6)\notin G$ and $(v_3,v_5)\notin G$; otherwise,
$N(v)\cap N(v_3)=\{u_2,u_5,u_6\}$ or $N(u_3)\cap
N(v_3)=\{v_1,u_5,v_5\}$, which is a contradiction. Therefore, $v_3$
contains a new neighbor, denoted by $v_6$. Since $v_3\in N(v_1)\cap
N(v_6)$, we have $(w,v_6)\in G$, since $w$ is the unique neighbor of
$v_1$ with degree less than $4$. Similarly, we get $(v_4,v_6)\in G$
by considering $N(v_2)\cap N(v_6)$. We further obtain that
$(w,v_5)\in G$ by considering $N(v_2)\cap N(v_5)$.

Up to now, $d(u_5)=d(u_6)=d(v_5)=d(v_6)=3$ and other vertices above
have degree $4$. It is known that $G[\{u_5,u_6,v_5,v_6\}]$ contains
no edges because of the triangle-free property of $G$. Suppose now
that $s$ is the forth neighbor of $u_5$. Considering $N(v)\cap
N(s)$, $N(u_3)\cap N(s)$ and $N(v_3)\cap N(s)$, respectively, we
derive that $(u_5,s)\in G, (u_6,s)\in G, (v_5,s)\in G$ and
$(v_6,s)\in G$. Now all vertices have degree $4$. It can be verified
that $G$ is the hypercube $Q_4$.

On the other hand, $(v_4,u_5)\in G$. It follows that $v_4 \in
N(v_2)\cap N(u_5)$. Then we have $(w,u_5)\in G$, since $w$ is the
unique neighbor of $v_2$ with degree less than $4$ other than $v_4$.
Note that $u_2\in N(v)\cap N(v_3)$, which forces that $u_6$ becomes
another common neighbor between $v$ and $v_3$, since $u_6$ is the
unique neighbor of $v$, whose degree is less than $4$. By a similar
discussion on $N(u_3)\cap N(v_3)$, we can deduce that $(v_3,v_5)\in
G$. Since $u_5\in N(u_3)\cap N(v_4)$, we get $(v_4,v_5)\in G$. We
further consider $N(u_4)\cap N(w)$ and obtain $(w,u_6)\in G$. Now
all vertices have degree $4$. It can be easily verified that $G$ is
the graph $G_3$ depicted in Figure \ref{Fig4}.

\textbf{Case 2.} $(a_1,a_2,a_3)=(3,1,1)$.

In this case, $v$ is adjacent to all vertices of $\{u_2,v_1,v_2\}$,
while $u_3$ and $u_4$ are adjacent to one of them, respectively.
Without loss of generality, $(u_2,u_3)\in G$. Then $(v_1,u_3)\notin
G$ and $(v_2,u_3)\notin G$ since $a_2=1$. It follows that
$(u_2,u_4)\in G$. If not, then either $(v_1,u_4)\in G$ or
$(v_2,u_4)\in G$, where the former possibility implies that
$N((u_1)\cap N(v_1)=\{u,v,u_4\})$ and the latter implies that
$N(u_1)\cap N(v_2)=\{u,v,u_4\}$, both of which contradict
Proposition \ref{4neigh}. Hence $(u_2,u_4)\in G$, $(v_1,u_4)\notin
G$ and $(v_2,u_4)\notin G$. Now let $v_3$ and $v_4$ be another two
neighbors of $v_1$. Observe that $v_1\in N(v)\cap N(v_3)$, which
forces $v_2$ to be another common vertex between $v$ and $v_3$. By a
similar discussion on $N(v)\cap N(v_4)$, we can derive $(v_2,v_4)\in
G$.

Now, $d(u)=d(v)=d(u_1)=d(u_2)=d(v_1)=d(v_2)=4$ and $d(u_3)=d(u_4)=d(v_3)=d(v_4)=2$.
We can divide our subsequent discussion into the following steps.
\begin{enumerate}[\indent Step $1.$]
\item If the induced subgraph $G[\{u_3,u_4,v_3,v_4\}]$ contains edges,
      then the edges can only be some of $(u_3,v_3)$, $(u_4,v_4)$, $(u_3,v_4)$
      and $(u_4,v_3)$, since $G$ is triangle-free. Without loss generality,
      assume $(u_3,v_4)\in G$. Then $u_3\in N(u_1)\cap N(v_4)$ and
      $v_4\in N(v_2)\cap N(u_3)$, which forces $(u_4,v_4)\in G$ and $(u_3,v_3)\in G$.
      We further consider $N(u_2)\cap N(v_3)$, and obtain $(u_4,v_3)\in G$.
      Consequently, each vertex above has degree $4$. It is easy to verify
      that $G$ is the graph $\mathcal{G}_1$ depicted in Figure \ref{Fig1}.
      Now the discussion stop;
\item If the induced subgraph $G[\{u_3,u_4,v_3,v_4\}]$ contains no edges,
      then there are another two neighbors of $u_3$, say $u_5$ and $u_6$.
      Considering $N(u_1)\cap N(u_5)$ and $N(u_1)\cap N(u_6)$, respectively, we have
      $(u_4,u_5)\in G$ and $(u_4,u_6)\in G$, since $u_4$ is the unique neighbor
      of $u_1$ whose degree is less than $4$.\\
      On the one hand, if $u_5$ or $u_6$ is adjacent to $v_3$ or $v_4$, then without loss generality we
      can suppose $(u_5,v_3)\in G$. Then $v_3\in N(v_1)\cap N(u_5)$,
      which implies $(u_5,v_4)\in G$. Notice that $u_5\in N(u_3)\cap N(v_3)$
      and $u_5\in N(u_3)\cap N(v_4)$. Then we deduce that $(v_3,u_6)\in G$
      and $(v_4,u_6)\in G$, since $u_6$ is the unique neighbor of $u_3$ whose
      degree is less than $4$. It can be verified that $G$ is the graph $\mathcal{H}_2$
      depicted in Figure \ref{Fig2}.\\
      On the other hand, both $u_5$ and $u_6$ are not adjacent to $v_3$ or $v_4$.
      Then $v_3$ has another two neighbors, denoted by $v_5$ and $v_6$.
      By a similar discussion on $N(v_2)\cap N(v_5)$ and $N(v_2)\cap N(v_6)$,
      respectively, we can obtain $(v_4,v_5)\in G$ and $(v_4,v_6)\in G$.
      Then continue the following step;
\item If the induced subgraph $G[\{u_5,u_6,v_5,v_6\}]$ contains edges,
      we can discuss this case similar to Step $1$. Consequently,
      we can obtain that $G$ is the graph $\mathcal{G}_2$ depicted in Figure \ref{Fig1}.
      The discussion stops; If the induced subgraph $G[\{u_5,u_6,v_5,v_6\}]$
      contains no edges, we can also continue the discussion according to Step $2$,
      until we get that $G$ is the graph $\mathcal{H}_3$ depicted in Figure
      \ref{Fig2}, or executing Step $3$ again. The discussion continues.
\end{enumerate}

It should be pointed out that the discussion will terminate by
illustrating that $G$ is either a graph in $\mathcal{G}_i$ or a
graph in $\mathcal{H}_j$, which are shown in Figure \ref{Fig1} and
Figure \ref{Fig2}, respectively.

\textbf{Case 3.} $(a_1,a_2,a_3)=(3,3,1)$.

This case means that $v$ and $u_3$ are adjacent to all vertices of
$\{u_2,v_1,v_2\}$, while $u_4$ is precisely adjacent to one of them.
Without loss generality, suppose $(u_2,u_4)\in G$. Then
$(v_1,u_4)\notin G$ and $(v_2,u_4)\notin G$. Consequently,
$|N(u_1)\cap N(v_1)|=|\{u,v,u_3\}|=3$, which contradicts Proposition
\ref{4neigh}. Therefore, this case could not happen.

\textbf{Case 4.} $(a_1,a_2,a_3)=(3,3,3)$.

Obviously, $v$, $u_3$ and $u_4$ are adjacent to all vertices of $\{u_2,v_1,v_2\}$.
It can be checked that all vertices have degree $4$, and hence $G$ is the
complete bipartite graph $K_{4,4}$, which is also the graph $\mathcal{H}_1$
depicted in Figure \ref{Fig2}.

To sum up the discussion above, $G$ is the hypercube $Q_4$ or the
graph $G_3$ or a graph in $\mathcal{G}_i$ or a graph in
$\mathcal{H}_j$. The proof is now complete.\qed

For convenience, we denote the set of all graphs presented above by
$\mathcal{F}$, which consists of $G_1$, $G_2$, $G_3$, $Q_4$, all
graphs in $\mathcal{G}_i$ and all graphs in $\mathcal{H}_j$.
Combining Theorem \ref{underlying1} with Theorem \ref{underlying2},
we conclude one of our main results as follows.

\begin{thm}\label{underlying}
Let $G^\sigma$ be a $4$-regular oriented graph with optimum skew
energy. Then the underlying graph $G$ is a graph in $\mathcal{F}$.
\end{thm}

Now the question naturally arises: whether there exists an orientation
for each graph of $\mathcal{F}$ such that the resultant oriented graph
attains optimum skew energy. The following results tell us that for
each graph of $\mathcal{F}$ such orientation indeed exists.

\begin{figure}[h,t,b,p]
\begin{center}
\scalebox{1}[1]{\includegraphics{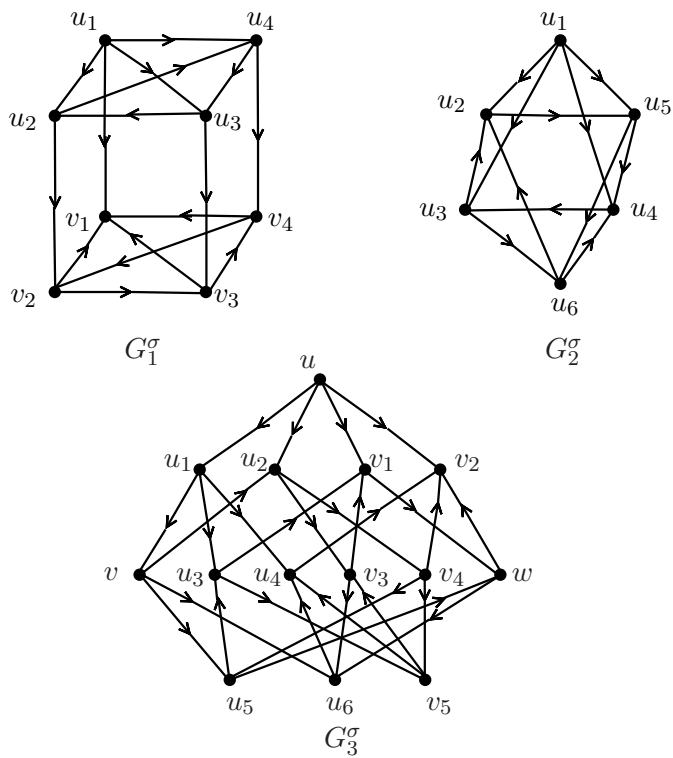}}
\end{center}
\caption{The optimum orientations for $G_1$, $G_2$ and $G_3$}\label{Fig5}
\end{figure}

\begin{thm}\label{orientation1}
Let $G_1^\sigma$, $G_2^\sigma$ and $G_3^\sigma$ be the oriented
graphs of $G_1$, $G_2$ and $G_3$, respectively, given in Figure
\ref{Fig5}. Then each of them has the optimum shew energy.
\end{thm}
\pf Let the rows of the skew-adjacency matrix $S(G_1^\sigma)$ correspond
successively the vertices $u_1$, $u_2$, $u_3$, $u_4$, $v_1$, $v_2$, $v_3$
and $v_4$. It follows that

\[
S(G_1^\sigma)=\begin{bmatrix}
0 & 1 & 1 & 1 & 1 & 0 & 0 & 0  \\-1 & 0 & -1 & 1 & 0 & 1 & 0 & 0\\
-1 & 1 & 0 & -1 & 0 & 0 & 1 & 0 \\-1 & -1 & 1 & 0 & 0 & 0 & 0 & 1\\
-1 & 0 & 0 & 0 & 0 & -1 & -1 & -1\\0 & -1 & 0 & 0 & 1 & 0 & 1 & -1\\
0 & 0 & -1 & 0 & 1 & -1 & 0 & 1  \\0 & 0 & 0 & -1 & 1 & 1 & -1 & 0\\
\end{bmatrix}
\]
Let the rows of the skew-adjacency matrix $S(G_2^\sigma)$ correspond
successively the vertices $u_1$, $u_2$, $u_3$, $u_4$, $u_5$ and $u_6$. Then
\[
S(G_2^\sigma)=\begin{bmatrix}
0 & 1 & 1 & 1 & 1 & 0 \\-1 & 0 & -1 & 0 & 1 & -1\\
-1 & 1 & 0 & -1 & 0 & 1\\-1 & 0 & 1 & 0 & -1 & -1\\
-1 & -1 & 0 & 1 & 0 & 1 \\0 & 1 & -1 & 1 & -1 & 0 \\
\end{bmatrix}
\]
Similarly, let the rows of the skew-adjacency matrix $S(G_3^\sigma)$
correspond successively the vertices $u$, $u_1$, $u_2$, $v_1$, $v_2$,
$v$, $u_3$, $u_4$,  $v_3$, $v_4$, $w$, $u_5$, $u_6$ and $v_5$. Then

\setcounter{MaxMatrixCols}{20}
\[
S(G_3^\sigma)=\begin{bmatrix}
0 & 1 & 1 & 1 & 1 & 0 & 0 & 0 & 0 & 0 & 0 & 0 & 0 & 0\\
-1 & 0 & 0 & 0 & 0 & 1 & 1 & 1 & 0 & 0 & 0 & 0 & 0 & 0\\
-1 & 0 & 0 & 0 & 0 & -1 & 0 & 0 & 1 & 1 & 0 & 0 & 0 & 0\\
-1 & 0 & 0 & 0 & 0 & 0 & -1 & 0 & -1 & 0 & 1 & 0 & 0 & 0\\
-1 & 0 & 0 & 0 & 0 & 0 & 0 & -1 & 0 & -1 & -1 & 0 & 0 & 0\\
0 & -1 & 1 & 0 & 0 & 0 & 0 & 0 & 0 & 0 & 0 & 1 & 1 & 0\\
0 & -1 & 0 & 1 & 0 & 0 & 0 & 0 & 0 & 0 & 0 & -1 & 0 & 1\\
0 & -1 & 0 & 0 & 1 & 0 & 0 & 0 & 0 & 0 & 0 & 0 & -1 & -1\\
0 & 0 & -1 & 1 & 0 & 0 & 0 & 0 & 0 & 0 & 0 & 0 & 1 & -1\\
0 & 0 & -1 & 0 & 1 & 0 & 0 & 0 & 0 & 0 & 0 & 1 & 0 & 1\\
0 & 0 & 0 & -1 & 1 & 0 & 0 & 0 & 0 & 0 & 0 & -1 & 1 & 0\\
0 & 0 & 0 & 0 & 0 & -1 & 1 & 0 & 0 & -1 & 1 & 0 & 0 & 0\\
0 & 0 & 0 & 0 & 0 & -1 & 0 & 1 & -1 & 0 & -1 & 0 & 0 & 0\\
0 & 0 & 0 & 0 & 0 & 0 & -1 & 1 & 1 & -1 & 0 & 0 & 0 & 0\\
\end{bmatrix}
\]
It is not difficult to check that
$S(G_1^\sigma)^TS(G_1^\sigma)=4I_8$,
$S(G_2^\sigma)^TS(G_2^\sigma)=4I_6$ and
$S(G_3^\sigma)^TS(G_3^\sigma)=4I_{14}$. We can also verify these
equalities by proving that different row vectors of each of
$S(G_1^\sigma)$, $S(G_2^\sigma)$ and $S(G_3^\sigma)$ are pairwise
orthogonal. The theorem is thus proved.\qed

We have known from \cite{Tian} that there exists an orientation
$\sigma$ of $Q_4$ such that the resultant oriented graph
$Q_4^\sigma$ has optimum skew energy. The following two algorithms
recursively describe optimum orientations of $\mathcal{G}_i$ and
$\mathcal{H}_j$, respectively.

\begin{figure}[h,t,b,p]
\begin{center}
\scalebox{1}[1]{\includegraphics{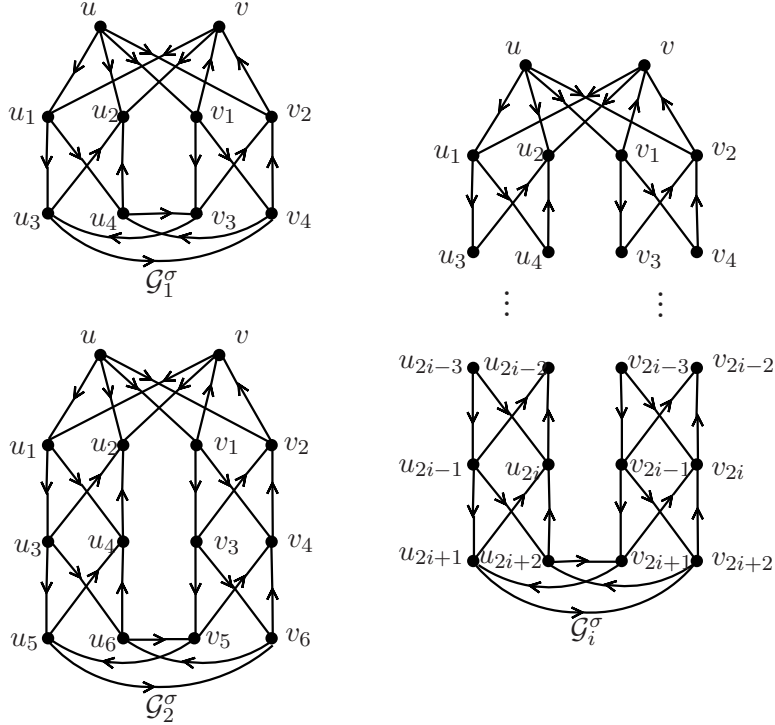}}
\end{center}
\caption{The optimum orientation for $\mathcal{G}_i$}\label{Fig6}
\end{figure}

\noindent\textbf{Algorithm 1.}
\begin{enumerate}[\indent Step $1.$]
\item Give $\mathcal{G}_1$ an orientation as shown in Figure \ref{Fig6}.
\item Assume that $\mathcal{G}_1,\mathcal{G}_2,\dots,\mathcal{G}_{t-1}$
      have been oriented into $\mathcal{G}_1^\sigma,\mathcal{G}_2^\sigma,\dots,\mathcal{G}_{t-1}^\sigma$.
      Then we orient $\mathcal{G}_t$ with the following method:
\begin{enumerate}[\indent(i)]
\item Keep the orientations of all edges in $E(\mathcal{G}_{t-1})\cap E(\mathcal{G}_t)$.
\item Give the remaining edges orientations such that $\langle u_{2t-1},u_{2t+1}\rangle$,
      $\langle u_{2t-1},u_{2t+2}\rangle$, $\langle u_{2t+1},u_{2t}\rangle$,
      $\langle u_{2t+2},u_{2t}\rangle$, $\langle v_{2t-1},v_{2t+1}\rangle$,
      $\langle v_{2t-1},v_{2t+2}\rangle$, $\langle v_{2t+1},v_{2t}\rangle$,
      $\langle v_{2t+2},v_{2t}\rangle$, $\langle u_{2t+1},v_{2t+2}\rangle$,
      $\langle v_{2t+2},u_{2t+2}\rangle$, $\langle u_{2t+2},v_{2t+1}\rangle$ and
      $\langle v_{2t+1},u_{2t+1}\rangle$ belong to $\Gamma(\mathcal{G}_t^\sigma)$.
\end{enumerate}
\item If $t=i$, stop; else take $t-1:=t$, return to Step 2.
\end{enumerate}

\begin{figure}[h,t,b,p]
\begin{center}
\scalebox{1}[1]{\includegraphics{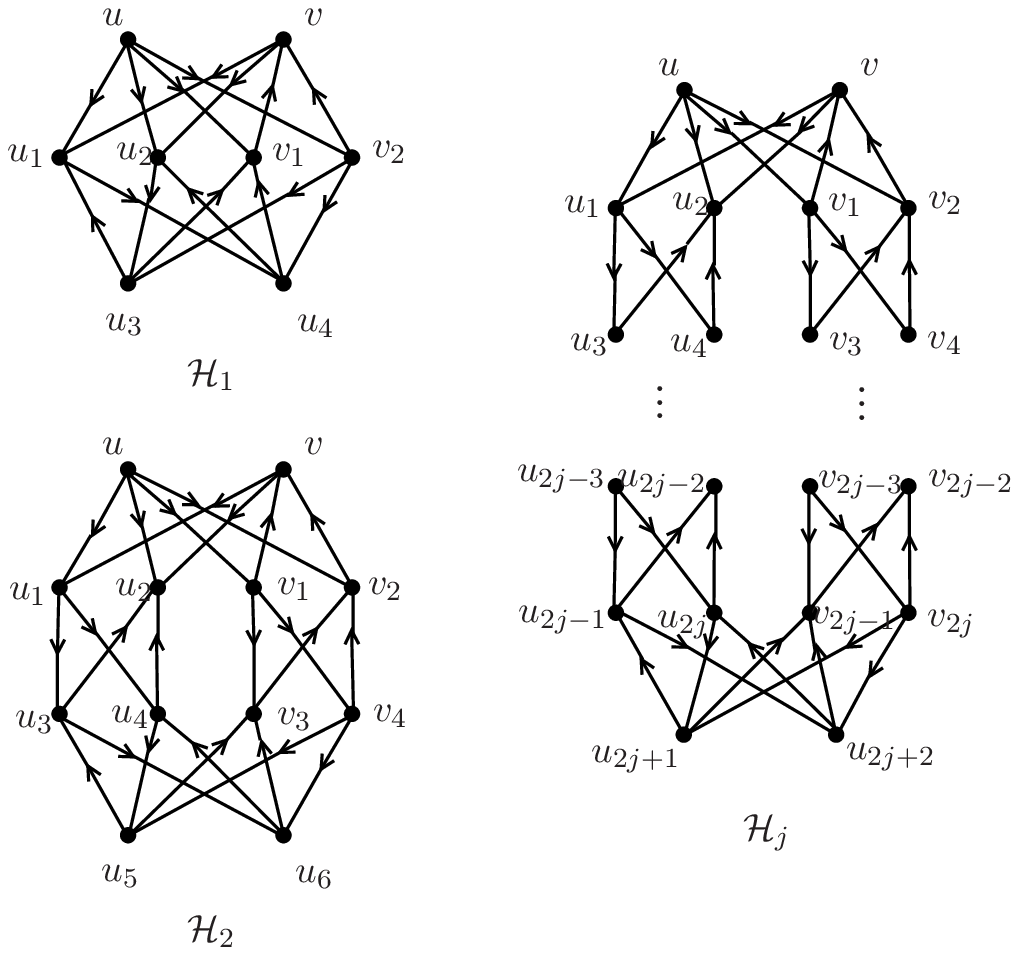}}
\end{center}
\caption{The optimum orientation for $\mathcal{H}_j$}\label{Fig7}
\end{figure}

\noindent\textbf{Algorithm 2.}
\begin{enumerate}[\indent Step $1.$]
\item Give $\mathcal{H}_1$ an orientation as shown in Figure \ref{Fig7}.
\item Assume that $\mathcal{H}_1,\mathcal{H}_2,\dots,\mathcal{H}_{t-1}$
      have been oriented into $\mathcal{H}_1^\sigma,\mathcal{H}_2^\sigma,\dots,\mathcal{H}_{t-1}^\sigma$.
      Then we orient $\mathcal{H}_t$ with the following method:
\begin{enumerate}[\indent(i)]
\item Keep the orientations of all edges in $E(\mathcal{H}_{t-1})$$\cap E(\mathcal{H}_t)$$\backslash$$ \{(u_{2t-3},u_{2t-1}),(u_{2t-3},u_{2t}),$\\$(u_{2t-2},u_{2t-1}),(u_{2t-2},u_{2t})\}$.
\item Give the remaining edges orientations such that $\langle u_{2t-3},u_{2t-1}\rangle$,
      $\langle u_{2t-3},u_{2t}\rangle$, $\langle u_{2t-1},u_{2t-2}\rangle$,
      $\langle u_{2t},u_{2t-2}\rangle$, $\langle v_{2t-3},v_{2t-1}\rangle$,
      $\langle v_{2t-3},v_{2t}\rangle$, $\langle v_{2t-1},v_{2t-2}\rangle$,
      $\langle v_{2t},v_{2t-2}\rangle$, $\langle u_{2t+1},u_{2t-1}\rangle$,
      $\langle u_{2t-1},u_{2t+2}\rangle$, $\langle u_{2t},u_{2t+1}\rangle$,
      $\langle u_{2t+2},u_{2t}\rangle$, $\langle u_{2t+1},v_{2t-1}\rangle$,
      $\langle u_{2t+2},v_{2t-1}\rangle$ $\langle v_{2t},u_{2t+1}\rangle$ and
      $\langle v_{2t},u_{2t+2}\rangle$ belong to $\Gamma(\mathcal{H}_t^\sigma)$.
\end{enumerate}
\item If $t=i$, stop; else take $t-1:=t$, return to Step 2.
\end{enumerate}

Next, we shall prove that $\mathcal{G}_i^\sigma$ and
$\mathcal{H}_j^\sigma$ derived from Algorithm 1 and Algorithm 2,
respectively, have optimum skew energy, that is, their
skew-adjacency matrices satisfy
$S(\mathcal{G}_i^\sigma)^TS(\mathcal{G}_i^\sigma)=4I$ and
$S(\mathcal{H}_j^\sigma)^TS(\mathcal{H}_j^\sigma)=4I$. In order to
illustrate clearly the skew-adjacency matrices
$S(\mathcal{G}_i^\sigma)$ and $S(\mathcal{H}_j^\sigma)$, we here
define some small matrix blocks.
\begin{alignat*}{2}
A&=\begin{bmatrix}
1 & 1 & 1 & 1\\1 & 1 & -1 & -1\\
\end{bmatrix}&\hspace{25pt}
B&=\begin{bmatrix}
1 & 1 & 0 & 0\\-1 & -1 & 0 & 0\\
0 & 0 & 1 & 1\\0 & 0 & -1 & -1\\
\end{bmatrix}\\
C&=\begin{bmatrix}
0 & 0 & -1 & 1\\0 & 0 & 1 & -1\\
1 & -1 & 0 & 0\\-1 & 1 & 0 & 0\\
\end{bmatrix}&\hspace{25pt}
D&=\begin{bmatrix}
-1 & 1\\ 1 & -1\\-1 & -1\\1 & 1\\
\end{bmatrix}
\end{alignat*}

\begin{thm}\label{orientation2}
Let $S(\mathcal{G}_i^\sigma)$ be the skew-adjacency matrix of $\mathcal{G}_i^\sigma$
obtained from Algorithm 1. Then $S(\mathcal{G}_i^\sigma)^TS(\mathcal{G}_i^\sigma)=4I$.
\end{thm}
\pf Let the rows of the skew-adjacency matrix $S(\mathcal{G}_i)$
correspond successively the vertices $u$, $v$, $u_1$, $u_2$, $v_1$,
$v_2$, $\dots$, $u_{2i+1}$, $u_{2i+2}$, $v_{2i+1}$, $v_{2i+2}$.
Then from Algorithm 1, $S(\mathcal{G}_i)$ can be written as the block
matrix for each positive integer $i$.

For $i=1$ and $i=2$,
\begin{alignat*}{2}
S(\mathcal{G}_1^\sigma)=\begin{bmatrix}
0 & A & 0\\
-A^T & 0 & B\\
0 & -B^T & C
\end{bmatrix} & \hspace{25pt}
S(\mathcal{G}_2^\sigma)=\begin{bmatrix}
0 & A & 0 & 0\\
-A^T & 0 & B &0\\
0 & -B^T & 0 & B\\
0 & 0 & -B^T & C
\end{bmatrix}
\end{alignat*}
By applying multiplication of block matrix, it is easy to compute
that
\begin{alignat*}{1}
S(\mathcal{G}_1^\sigma)^TS(\mathcal{G}_1^\sigma)&=\begin{bmatrix}
AA^T & 0 & -AB\\
0 & A^TA+BB^T & -BC\\
-B^TA^T & -C^TB^T & B^TB+C^TC
\end{bmatrix}\\
S(\mathcal{G}_2^\sigma)^TS(\mathcal{G}_2^\sigma)&=\begin{bmatrix}
AA^T & 0 & -AB & 0\\
0 & A^TA+BB^T & 0 & -B^2\\
-B^TA^T & 0 & B^TB+BB^T & -BC\\
0 & -(B^T)^2 & -C^TB^T & B^TB+C^TC
\end{bmatrix}
\end{alignat*}

In order to prove
$S(\mathcal{G}_1^\sigma)^TS(\mathcal{G}_1^\sigma)=4I$ and
$S(\mathcal{G}_2^\sigma)^TS(\mathcal{G}_2^\sigma)=4I$, it suffices
to prove that the following equalities meanwhile hold.
\begin{equation}\label{block1}
\begin{split}
&AA^T=4I_2,\,\, A^TA+BB^T=4I_4,\,\, B^TB+BB^T=4I_4,\\
& B^TB+C^TC=4I_4,\,\, AB=0,\,\, B^2=0,\,\, BC=0.
\end{split}
\end{equation}
By the definitions of $A$, $B$ and $C$, it is easy to verify that
all equalities \ref{block1} indeed hold. In fact, these equalities
can further guarantee
$S(\mathcal{G}_i^\sigma)^TS(\mathcal{G}_i^\sigma)=4I$, because
$S(\mathcal{G}_i)$ can be formulated as
\begin{equation*}
S(\mathcal{G}_i^\sigma)=\left[\begin{array}{ccccccc}
0 & A & 0 & 0 & \cdots & 0 & 0\\
-A^T & 0 & B & 0 & \cdots & 0 & 0\\
0 & -B^T & 0 & B & \cdots & 0 & 0\\
0 & 0 & -B^T & 0 & \cdots & 0 & 0\\
\vdots & \vdots & \vdots & \vdots & \vdots & \vdots & \vdots\\
0 & 0 & 0 & 0 & \cdots & 0 & B\\
0 & 0 & 0 & 0 & \cdots & -B^T & C\\
\end{array}\right].
\end{equation*}

It should be pointed out that if one only considers $\mathcal{G}_1$,
then it is enough to check that parts of the equalities hold. The
proof is now complete.\qed

\begin{thm}\label{orientation3}
Let $S(\mathcal{H}_j^\sigma)$ be the skew-adjacency matrix of $\mathcal{H}_j^\sigma$
obtained from Algorithm 2. Then $S(\mathcal{H}_j^\sigma)^TS(\mathcal{H}_j^\sigma)=4I$.
\end{thm}
\pf Similar to the proof of Theorem \ref{orientation2}, let the rows of
the skew-adjacency matrix $S(\mathcal{H}_j)$ correspond successively
the vertices $u$, $v$, $u_1$, $u_2$, $v_1$, $v_2$, $\dots$, $u_{2i-1}$,
$u_{2i}$, $v_{2i-1}$, $v_{2i}$, $u_{2i+1}$ and $u_{2i+2}$.
\begin{alignat*}{2}
S(\mathcal{H}_1^\sigma)=\begin{bmatrix}
0 & A & 0\\
-A^T & 0 & D\\
0 & -D^T & 0
\end{bmatrix} & \hspace{25pt}
S(\mathcal{H}_2^\sigma)=\begin{bmatrix}
0 & A & 0 & 0\\
-A^T & 0 & B &0\\
0 & -B^T & 0 & D\\
0 & 0 & -D^T & 0
\end{bmatrix}
\end{alignat*}
\begin{equation*}
S(\mathcal{H}_j^\sigma)=\left[\begin{array}{cccccccc}
0 & A & 0 & 0 & \cdots & 0 & 0 & 0\\
-A^T & 0 & B & 0 & \cdots & 0 & 0 & 0\\
0 & -B^T & 0 & B & \cdots & 0 & 0 & 0\\
0 & 0 & -B^T & 0 & \cdots & 0 & 0 & 0\\
\vdots & \vdots & \vdots & \vdots & \vdots & \vdots & \vdots & \vdots\\
0 & 0 & 0 & 0 & \cdots & 0 & B & 0\\
0 & 0 & 0 & 0 & \cdots & -B^T & 0 & D\\
0 & 0 & 0 & 0 & \cdots & 0 & -D^T & 0\\
\end{array}\right]
\end{equation*}

We can verify that $S(\mathcal{H}_1)^TS(\mathcal{H}_1)=4I$ if and only if the equalities below hold,
\begin{equation}\label{block2}
AA^T=4I_2,\,\, A^TA+DD^T=4I_4,\,\,D^TD=4I_2,\,\, AD=0.
\end{equation}
while $S(\mathcal{H}_2)^TS(\mathcal{H}_2)=4I$ if and only if the following equalities hold,
\begin{equation}\label{block3}
\begin{split}
&AA^T=4I_2,\,\, A^TA+BB^T=4I_4,\,\,B^TB+D^TD=4I_4,\\
&D^TD=4I_2,\,\, AB=0,\,\,BD=0.
\end{split}
\end{equation}
For $j\geq 3$, combining equalities in (\ref{block3}) with equalities
$B^TB+BB^T=4I_4$ and $B^2=0$, it is enough to ensure that the equality
$S(\mathcal{H}_j^\sigma)^TS(\mathcal{H}_j^\sigma)=4I$ holds.

By the definitions of $A,B$ and $D$, it can be directly checked that
the all equalities above indeed hold. This completes the proof.\qed

We can summarize all results above as the following theorem.
\begin{thm}
Let $G$ be a $4$-regular graph. Then $G$ has an optimum orientation
if and only if $G$ is a graph of $\mathcal{F}$.
\end{thm}

\textbf{Remark 1.} For arbitrary matrices $A'$, $B'$, $C'$ and $D'$
with entries $0$, $1$ and $-1$, if they have the same orders and the
same number of $0'$s with $A$, $B$, $C$ and $D$, respectively, and
meanwhile they satisfy all the equalities of Theorem
\ref{orientation2} and Theorem \ref{orientation3}, then we can
substitute $A$, $B$, $C$ and $D$, respectively by $A'$, $B'$, $C'$
and $D'$ in the skew-adjacency matrices $S(\mathcal{G}_i)$ and
$S(\mathcal{H}_j)$, and the corresponding oriented graphs still have
optimum skew energy.

\textbf{Remark 2.} The proofs of Theorem \ref{orientation2} and
Theorem \ref{orientation3} are based on matrix computations by
proving that the skew-adjacency matrix $S$ satisfies $S^TS=nI$.
Besides, we can apply Proposition \ref{2-walk} to prove that for any
two distinct vertices $u$ and $v$, the number of all positive walks
equals that of all negative walks from $u$ to $v$ with length $2$.


\begin{thebibliography}{20}

\bibitem{ABC}C. Adiga, R. Balakrishnan, W. So, The shew energy of a digraph, Linear Algebra Appl. 432(2010), 1825--1835.

\bibitem{CLL}X. Chen, X. Li, H. Lian, The skew energy of random oriented graphs, Linear Algebra Appl. 438(2013), 4547--4556.

\bibitem{GX}S. Gong, G. Xu, $3$-Regular digraphs with optimum skew energy, Linear Algebra Appl. 436(2012), 465--471.

\bibitem{G}I. Gutman, The energy of a graph, Ber. Math. Statist. Sekt. Forschungsz. Graz, 103(1978), 1--22.

\bibitem{GLZ}I. Gutman, X. Li, J. Zhang, Graph Energy, in: M. Dehmer, F. Emmert-Streib (Eds.), Analysis of Complex Network:
From Biology to Linguistics, Wiley-VCH Verlag, Weinheim, 2009, 145--174.

\bibitem{HL}Y. Hou, T. Lei, Characteristic polynomials of skew-adjacency matrices of oriented graphs, Electron. J. Combin.
18(2011), 156--167.

\bibitem{HSZ}Y. Hou, X. Shen, C. Zhang, Oriented unicyclic graphs with extremal skew energy,
   Available at http://arxiv.org/abs/1108.6229.

\bibitem{LSG}X. Li, Y. Shi, I. Gutman, Graph Energy, Springer, New York, 2012.

\bibitem{Shader} B. Shader, W. So, Skew spectra of oriented graphs, Electron. J. Combin. 16(2009), \#N32.

\bibitem{HSZ2}X. Shen, Y. Hou, C. Zhang, Bicyclic digraphs with extremal skew energy, Electron. J. Linear Algebra 23(2012), 340--355.

\bibitem{Tian}G. Tian, On the skew energy of orientations of hypercubes, Linear Algebra Appl. 435(2011), 2140--2149.

\end{thebibliography}
\end{document}